\documentclass[12pt]{article}

\usepackage{graphicx, xcolor}
\usepackage{epsfig}
\usepackage{amssymb}
\usepackage[T1]{fontenc}

\def\i{{\mathbf i}}

\def\e{{\mathbf e}}

 \def\R{\mathbb{R}}

\def\C{\mathbb{C}}
\def\N{\mathbb{N}}

\newcommand{\proof}{\bf {Proof:} \rm}

\newtheorem{theorem}{Theorem}[section]
\newtheorem{remark}{Remark}[section]
\newtheorem{lemma}{Lemma}[section]
\newtheorem{proposition}{Proposition}[section]

\begin{document}

\title{A Generalization of Fueter's Theorem in Dunkl-Clifford Analysis}

\author{Shanshan Li\thanks{School of Computer Science and Technology, Southwest University for
Nationalities, Chengdu, 610041, P. R. China. E-Mail:
shanshanli@swun.cn} \  and Minggang Fei\thanks{School of
Mathematical Sciences, University of Electronic Science and
Technology of China,
 Chengdu, 610054, P. R. China. E-Mail: fei@uestc.edu.cn} \thanks{Corresponding author}}
\date{}

\maketitle

\begin{abstract}In this paper we first offer an alternative approach to extend the original Fueter's Theorem in Dunkl-Clifford analysis to a version of the higher order case. Then this result is used to prove a generlized version of Fueter's Theorem with an extra homogeneous
Dunkl-monogenic polynomial $P_n(x_0,\underline{x})$ in $\R_1^d$
instead of the classical factor $P_n(\underline{x})$ in $\R^d$.
\end{abstract}

{\bf MSC 2000}: 30G35, 31A05, 33C80

{\bf Key words}: Reflection group, Dunkl-Dirac operator,
Dunkl-monogenic function, Fueter's Theorem

\section{Introduction}

To explain the idea we start from the very basic facts. Let $f$ be a
holomorphic function in an open set of the upper half complex plane
and write
\begin{eqnarray*}
f(z)=u(s,t)+\i v(s,t), \qquad z=s+\i t,
\end{eqnarray*}
then the Fueter's Theorem in Dunkl-Clifford analysis(\cite{FCK3})
asserts that in the corresponding region there holds
\begin{eqnarray*}
D\Delta^{\gamma_{\kappa}+\frac{d-1}{2}}\left(u(x_0,|\underline{x}|)+\frac{\underline{x}}{|\underline{x}|}v(x_0,|\underline{x}|)\right)=0
\end{eqnarray*}
whenever $\gamma_{\kappa}+(d-1)/2$ is any positive integer. Here $\underline{x}\in\R^d$, $D=T_{x_0}+\underline{D}$, $\underline{D}=\sum_{i=1}^d\e_iT_{x_i}$ and $\Delta=\sum_{i=0}^dT^2_{x_i}$, where $T_{x_i}(i=0,1,\cdots,d)$ are Dunkl operators corresponding to a finite reflection group $W$ which leaves $x_0$-axis invariant.

Furthermore, the authors in \cite{FCK3} generalized the above result as follows:

If, under the same assumptions of the above theorem, $P_n(\underline{x})$ is a homogeneous Dunkl-monogenic function of degree $n$ with respect to the Dunkl-Dirac operator $\underline{D}$, i.e. $\underline{D}P_n(\underline{x})=0$, then
\begin{eqnarray*}
D\Delta^{\gamma_{\kappa}+n+\frac{d-1}{2}}\left(\left(u(x_0,|\underline{x}|)+\frac{\underline{x}}{|\underline{x}|}v(x_0,|\underline{x}|)\right)P_n(\underline{x})\right)=0
\end{eqnarray*}
whenever $\gamma_{\kappa}+n+(d-1)/2$ is a positive integer.

The goal of this paper is to prove that the above result is still
valid if we replace $P_n(\underline{x})$ by a homogeneous monogenic
polynomial $P_n(x_0,\underline{x})$ of degree $n$ in $\R^{d}_1$,
which is a counterpart of generalized Fueter's
Theorem(\cite{PS2010}) in classical Clifford analysis. Motivated by
\cite{PS2010}, we will first prove an extension of Fueter's Theorem
in Dunkl-Clifford analysis which uses complex-valued functions
satisfying the following equation
\begin{eqnarray}
\partial_{\bar{z}}\Delta_z^mf(z)=0, \ m\in\N_0,  \label{Eqn:1}
\end{eqnarray}
as initial function instead of the usual holomorphic functions,
where $z=t+\i s$, $\partial_{\bar{z}}$ and $\Delta_z$ denote,
respectively, the classical Cauchy-Riemann operator
\begin{eqnarray*}
\partial_{\bar{z}}=\frac{1}{2}(\partial_t+\i\partial_s)
\end{eqnarray*}
and Laplace operator in two dimensions
\begin{eqnarray*}
\Delta_z=\partial^2_t+\partial^2_s.
\end{eqnarray*}

The current paper is organized as follows. Section 2 contains some
basic facts about Dunkl-Clifford analysis and a characterization of
Dunkl-Dirac operator in spherical coordinates. A higher order
version of original Fueter's Theorem in Dunkl case is considered in
section 3. In the last section we prove a version of generalized
Fueter's Theorem with an extra Dunkl-monogenic factor
$P_n(x_0,\underline{x})$.

\section{Preliminaries and Dunkl-Clifford Analysis}

Let $\e_1, \cdots, \e_d$ be an orthonormal basis of $\R^d$
satisfying the anti-commutation relationship
$\e_i\e_j+\e_j\e_i=-2\delta_{ij}$, where $\delta_{ij}$ is the
Kronecker symbol. We define the universal real-valued Clifford
algebra $\R_{0,d}$(\cite{BDS},\cite{DSS}) as the $2^d$-dimensional
associative algebra with basis given by $\e_0=1$ and
$\e_A=\e_{l_1}\cdots \e_{l_n}$, where $A=\{l_1,\cdots,
l_n\}\subset\{1,\cdots,d\}$, for $1\leq l_1<\cdots<l_n\leq d$.
Hence, each element $x\in\R_{0,d}$ will be represented by
$x=\sum_Ax_A\e_A$, $x_A\in\R$.

In what follows, $sc[x]=x_0$ will denote the scalar part of
$x\in\R_{0,d}$, while an element $x=(x_0, x_1,\cdots,x_d)$ of
$\R^d_1$ will be identified with $x=x_0+\underline{x}$,
$\underline{x}=\sum_{i=1}^dx_i\e_i$. Also, we need the
anti-involution $\bar{\cdot}$ defined by $\bar{\e}_i=-\e_i$, and
$\overline{\e_i\e_j}=\bar{\e}_j\bar{\e}_i$. An important property of
algebra $\R_{0,d}$ is that each non-zero vector $x$ in $\R^d$(or in
$\R^d_1$) has a multiplicative inverse given by
$\frac{\bar{x}}{|x|^2}$. A $\R_{0,d}-$valued function $f$ over
$\Omega\subset\R_1^d$ has a representation $f=\sum_A\e_Af_A$ with
component $f_A: \Omega\rightarrow\R$.

The reflection $\sigma_\alpha x$ of a given vector $x\in\R^d_1$  on
the hyperplane $H_{\alpha}$ orthogonal to $\alpha\not=0$ is given by
$$
\sigma_\alpha x:=x-2\frac{\langle\alpha,x\rangle}{|\alpha|^2}\alpha.
$$
A finite set $R\subset\R^d_1\backslash\{0\}$ is called a root system
if $R\bigcap\R^d_1\cdot\alpha=\{\alpha,-\alpha\}$ and
$\sigma_{\alpha}R=R$ for all $\alpha\in R$. For a given root system
$R$ the reflections $\sigma_{\alpha}$, $\alpha\in R$, generate a
finite group $W\subset O(d)$, called the finite reflection group (or
Coxeter group) associated with $R$. All reflections in $W$
correspond to suitable pairs of roots. For a given
$\beta\in\R^d\backslash\bigcup_{\alpha\in R}H_{\alpha}$, we fix the
positive subsystem $R_{+}=\{\alpha\in
R|\langle\alpha,\beta\rangle>0\}$, i.e. for each $\alpha\in R$
either $\alpha\in R_{+}$ or $-\alpha\in R_{+}$. Sometimes we will
only consider reflections which only act in $\R^d$. In this case we
denote $\alpha$ or $\beta$ by $\underline{\alpha}$ or
$\underline{\beta}$. A function $\kappa: R\rightarrow\R^+$ on a root
system $R$ is called a multiplicity function if it is invariant
under the action of the associated reflection group $W$. For
abbreviation, we introduce the index
$\gamma_{\kappa}=\sum_{\alpha\in R_{+}}\kappa(\alpha)$ and the
Dunkl-dimension $\mu=2\gamma_{\kappa}+d$.

For each fixed positive subsystem $R_+$ and multiplicity function
$\kappa$ we have, as invariant operators, the
differential-difference operators (also called Dunkl operators)(\cite{Dunkl89},\cite{DX01}):
\begin{eqnarray}
T_{x_i}f(x)=\frac{\partial}{\partial x_i}f(x)+\sum_{\alpha\in
R_{+}}\kappa(\alpha)\frac{f(x)-f(\sigma_{\alpha}x)}{\langle\alpha,x\rangle}\alpha_i,
\qquad i=0, 1,  \cdots, d,
\end{eqnarray}
for $f\in C^1(\R^d_1)$. In the case $\kappa=0$, the $T_{x_i}$, $i=0,
1, \cdots, d$, reduce to the corresponding partial derivatives.
This also gives us the justification to think of these
differential-difference operators as the equivalent of partial
derivatives in the case of finite reflection groups. In this paper,
we will assume throughout that $\kappa\geq 0$ and
$\gamma_{\kappa}>0$. More importantly, these operators mutually
commute; that is, $T_{x_i}T_{x_j}=T_{x_j}T_{x_i}$. This property
allows us to define a Dunkl-Dirac operator in $\R^d$ for the
corresponding reflection group $W$ given by(\cite{CKR})
\begin{eqnarray}
\underline{D}f=\sum_{i=1}^d\e_iT_{x_i}f.
\end{eqnarray}
So, the Dunkl Laplacian $\underline{\Delta}$ in $\R^d$ is defined by
$\underline{\Delta}=-\underline{D}^2=\sum_{i=1}^dT_{x_i}^2$.

We now introduce the Dunkl-Cauchy-Riemann operator in $\R^d_1$
\begin{eqnarray*}
D=T_{x_0}+\underline{D},
\end{eqnarray*}
and Dunkl Laplacian in $\R^d_1$
\begin{eqnarray*}
\Delta=T_{x_0}^2+\underline{\Delta}.
\end{eqnarray*}
In this paper we will assume that our group $W$ will leave the
$x_0$-axis invariant. Since in this case we have
$T_{x_0}=\partial_{x_0}$ the Dunkl-Cauchy-Riemann operator and
Dunkl-Laplacian in $\R^d_1$ can also be written by
\begin{eqnarray}
D=\partial_{x_0}+\underline{D},
\end{eqnarray}
and
\begin{eqnarray}
\Delta=\partial^2_{x_0}+\underline{\Delta}.
\end{eqnarray}

Functions belonging to the kernel of Dunkl-Dirac operator
$\underline{D}$ or Dunkl-Cauchy-Riemann operator $D$ will be called
Dunkl-monogenic functions. As usual, functions belonging to be the
kernel of Dunkl Laplacian will be called Dunkl-harmonic functions.

From \cite{FCK3} we have the following
representation of the Dunkl-Dirac operator in spherical coordinates.
\begin{proposition} In spherical coordinates, i.e. $r=|\underline{x}|$ and $\underline{\omega}=\frac{\underline{x}}{|\underline{x}|}$ for $\underline{x}\in\R^d$, the Dunkl-Dirac operator has the form:
\begin{eqnarray}
\underline{D}=\underline{\omega}(\partial_r+\frac{1}{r}\Gamma_{\underline{\omega}}),
\end{eqnarray}
with
\begin{eqnarray*}
\Gamma_{\underline{\omega}}=\gamma_{\kappa}+\Phi_{\underline{\omega}}+\Psi,
\end{eqnarray*}
where
\begin{eqnarray*}
\Phi_{\underline{\omega}}=-\sum_{i<j}\e_i\e_j(x_i\partial_{x_j}-x_j\partial_{x_i}),
\end{eqnarray*}
and
\begin{eqnarray*}
\Psi f(\underline{x})=-\sum_{i<j}\e_i\e_j\sum_{\underline{\alpha}\in
R^+}\kappa(\underline{\alpha})\frac{f(\underline{x})-f(\sigma_{\underline{\alpha}}\underline{x})}{\langle\underline{\alpha},\underline{x}\rangle}(x_i\alpha_j-x_j\alpha_i)-\sum_{\underline{\alpha}\in
R^+}\kappa(\underline{\alpha})f(\sigma_{\underline{\alpha}}\underline{x}),
\end{eqnarray*}
for any $f\in C^1(\R^d)$.
\end{proposition}

\begin{remark} The operator $\Phi_{\underline{\omega}}$ in above proposition
corresponds to the classical spherical vector derivatives (the
classic Gamma operator) while the additional operator $\Psi$ and
constant $\gamma_{\kappa}$ derive from the difference part of Dunkl
operators.
\end{remark}

Furthermore, $\Gamma_{\underline{\omega}}$ is a first differential-difference
operator which satisfies the following properties(\cite{FCK3}):

\begin{proposition} Assume that $f(x)=f(|\underline{x}|)$ is a
radial function and $P_n(\underline{x})$ is a homogeneous Dunkl-monogenic
function of degree $n$, $n\in\R$, then
\begin{eqnarray*}
& &{\rm (i)} \ \Gamma_{\underline{\omega}} f(r)=0,\\
& &{\rm (ii)} \ \Gamma_{\underline{\omega}}P_n(\underline{\omega})=-nP_n(\underline{\omega}),\\
& &{\rm (iii)} \ \Gamma_{\underline{\omega}}(\underline{\omega}
P_n(\underline{\omega}))=(\mu+n-1)\underline{\omega}P_n(\underline{\omega}).
\end{eqnarray*}
\end{proposition}

Henceforward, we denote by ${\cal M}(n)$ the space of all
homogeneous Dunkl-monogenic polynomials of degree $n\in\N$. We then
immediately have
\begin{lemma}
Let $k\in\N$ and $P_n\in{\cal M}(n)$. Then there has
\begin{eqnarray*}
\underline{D}(\underline{x}^kP_{n}(\underline{x}))=\left\{
       \begin{array}{ll}
        -kx^{k-1}P_{n}(x), \ \  k \ even,\\
        \ \\
        -(k+\mu+2n-1)x^{k-1}P_{n}(x), \ \  k \ odd.
        \end{array}
        \right.
\end{eqnarray*}
\end{lemma}

\section{A higher order version of the original Fueter's Theorem}

In this section we prove an extension of the original Fueter's
Theorem in Dunkl case starting from a complex-valued function
satisfying equation (\ref{Eqn:1}) instead of the usual holomorphic
function. To this end, we start with the following lemma from
~\cite{PQS}, ~\cite{PP} or~\cite{PS2010} which we only state the
special case that we will use in this paper.

\begin{lemma}
Suppose that $f(x_0,r)$ is a scalar-valued infinitely differentiable
functions in $\R^2$ and that $D_r$ and $D^r$ are differential
operators defined by $D_r(0)\{f\}=D^r(0)\{f\}=f$ and
\begin{eqnarray*}
& &D_r(m)\{f\}=\left(\frac{1}{r}\partial_{r}\right)^m\{f\},\\
& &D^{r}(m)\{f\}=\partial_{r}\left(\frac{D^{r}(m-1)\{f\}}{r}\right)
\end{eqnarray*}
for integer $m\geq 1$. Then one has
\begin{eqnarray*}
& &{\rm (i)} \ \partial^2_{r}D_{r}(m)\{f\}=D_{r}(m)\{\partial^2_{r}f\}-2mD_{r}(m+1)\{f\},\\
& &{\rm (ii)} \ \partial^2_{r}D^{r}(m)\{f\}=D^{r}(m)\{\partial^2_{r}f\}-2mD^{r}(m+1)\{f\},\\
& &{\rm (iii)} \
D^{r}(m)\{\partial_{r}f\}=\partial_{r}D_{r}(m)\{f\},\\
& &{\rm (iv)} \
D_{r}(m)\{\partial_{r}f\}-\partial_{r}D^{r}(m)\{f\}=2m/rD^{r}(m)\{f\}.
\end{eqnarray*}
\end{lemma}


Furthermore, we need the following lemma which shows that the
iterated Dunkl-Laplacian $\Delta^m$, for any positive integer $m$,
keeps functions of the form
$(f(x_0,|\underline{x}|)+\frac{\underline{x}}{|\underline{x}|}g(x_0,|\underline{x}|))P_n(\underline{x})$
invariant whenever $f$ and $g$ are scalar-valued functions in $\R^2$
and $P_n(\underline{x})$ be a homogeneous Dunkl-monogenic function
of degree $n$ in $\R^d$.

\begin{lemma} Let $f(x_0,r)$ be a scalar-valued infinitely differentiable function
in an open set of $\R^2_+=\{(t,s)\in\R^2:s>0\}$. Then for $m\in\N_0$
we have
\begin{eqnarray*}
\Delta^m(f(x_0,|\underline{x}|)P_n(\underline{x}))=\left(\sum_{j=0}^md_{n,\mu}(j)\left(\begin{array}
{cc} m \\ j
\end{array}\right)D_r(j)\{\Delta_z^{m-j}f(x_0,r)\}\right)P_n(\underline{x})
\end{eqnarray*}
and
\begin{eqnarray*}
\Delta^m(f(x_0,|\underline{x}|)\frac{\underline{x}}{|\underline{x}|}P_n(\underline{x}))=\left(\sum_{j=0}^md_{n,\mu}(j)\left(\begin{array}
{cc} m \\ j
\end{array}\right)D^r(j)\{\Delta_z^{m-j}f(x_0,r)\}\right)\frac{\underline{x}}{|\underline{x}|}P_n(\underline{x}),
\end{eqnarray*}
where
\begin{eqnarray*}
& &d_{n,\mu}(0)=1,\\
& &d_{n,\mu}(j)=(2n+\mu-1)(2n+\mu-3)\cdots(2n+\mu-(2j-1)).
\end{eqnarray*}
\end{lemma}

\proof We will prove this lemma by induction. Let
$\underline{\omega}=\frac{\underline{x}}{|\underline{x}|}$. When
$m=1$, we need to show that the following identities hold
\begin{eqnarray}
\Delta(fP_n)=\left(\Delta_zf+(2n+\mu-1)D_r(1)\{f\}\right)P_n
\end{eqnarray}
and
\begin{eqnarray}
\Delta(f\underline{\omega}P_n)=(\Delta_zf+(2n+\mu-1)D^r(1)\{f\})\underline{\omega}P_n.
\end{eqnarray}

To prove (7), we start from
 \begin{eqnarray*}
\Delta_h=\partial^2_{x_0}-\underline{D}_h\underline{D}_h, \qquad
\underline{D}_h=\underline{\omega}(\partial_r+\frac{1}{r}\Gamma_{\underline{\omega}}).
\end{eqnarray*}
Then using Proposition 2.2 we get
\begin{eqnarray*}
\underline{D}(fP_n)&=&\underline{\omega}(\partial_r+\frac{1}{r}\Gamma_{\underline{\omega}})(fr^nP_n(\underline{\omega}))\\
                     &=&(\partial_rf)r^n\underline{\omega}P_n(\underline{\omega})
\end{eqnarray*}
and
\begin{eqnarray*}
\underline{D}\underline{D}(fP_n)&=&\underline{\omega}(\partial_r+\frac{1}{r}\Gamma_{\underline{\omega}})\left((\partial_rf)r^n\underline{\omega}P_n(\underline{\omega})\right)\\
                                    &=&-(\partial^2_rf)r^nP_n(\underline{\omega})-n(\partial_rf)r^{n-1}P_n(\underline{\omega})\\
                                    & &-(\mu+n-1)(\partial_rf)r^{n-1}P_n(\underline{\omega})\\
                                    &=&-\left(\partial^2_rf+\frac{2n+\mu-1}{r}(\partial_rf)\right)r^nP_n(\underline{\omega}).
\end{eqnarray*}
Therefore, we have
\begin{eqnarray*}
\Delta(fP_n)&=&\left(\partial^2_{x_0}f+\partial^2_rf+\frac{2n+\mu-1}{r}(\partial_rf)\right)r^nP_n(\underline{\omega})\\
              &=&\left(\Delta_zf+(2n+\mu-1)D_r(1)\{f\}\right)P_n.
\end{eqnarray*}

To prove (8), again applying Proposition 2.2 we obtain
\begin{eqnarray*}
\underline{D}(f\underline{\omega}P_n)&=&\underline{\omega}(\partial_r+\frac{1}{r}\Gamma_{\underline{\omega}})(fr^n\underline{\omega}P_n(\underline{\omega}))\\
                                       &=&\underline{\omega}((\partial_rf)r^n\underline{\omega}P_n(\underline{\omega})+nfr^{n-1}\underline{\omega}P_n(\underline{\omega})+fr^{n-1}\Gamma_{\underline{\omega}}(\underline{\omega}P_n(\underline{\omega})))\\
                                       &=&-(\partial_rf)r^nP_n(\underline{\omega})-(2n+\mu-1)fr^{n-1}P_n(\underline{\omega})
\end{eqnarray*}
and
\begin{eqnarray*}
& &\underline{D}\underline{D}(f\underline{\omega}P_n)\\
&=&-\underline{\omega}(\partial_r+\frac{1}{r}\Gamma_{\underline{\omega}})((\partial_rf)r^nP_n(\underline{\omega})+(2\gamma_{\kappa}+2n+d-1)fr^{n-1}P_n(\underline{\omega}))\\
                                                      &=&-\underline{\omega}(((\partial^2_rf)r^n+n(\partial_rf)r^{n-1})P_n(\underline{\omega})+(\partial_rf)r^{n-1}\Gamma_{\underline{\omega}}(P_n(\underline{\omega}))\\
                                                      & &\ \ +(2n+\mu-1)((\partial_rf)r^{n-1}+(n-1)fr^{n-2})P_n(\underline{\omega})\\
                                                      & &\ \ +(2n+\mu-1)fr^{n-2}\Gamma_{\underline{\omega}}(P_n(\underline{\omega})))\\
                                                      &=&-\left(\partial^2_rf+\frac{n}{r}(\partial_rf)-\frac{n}{r}(\partial_rf)+\frac{2n+\mu-1}{r}(\partial_rf)\right.\\
                                                      &
                                                      &+\left.\frac{(2n+\mu-1)(n-1)}{r^2}f-\frac{(2n+\mu-1)n}{r^2}f\right)\underline{\omega}r^nP_n(\underline{\omega})\\
                                                      &=&-\left(\partial^2_rf+(2n+\mu-1)(\frac{\partial_rf}{r}-\frac{f}{r^2})\right)\underline{\omega}P_n.
\end{eqnarray*}
This leads to
\begin{eqnarray*}
\Delta(f\underline{\omega}P_n)&=&\left(\partial^2_{x_0}f+\partial^2_rf+(2n+\mu-1)(\frac{\partial_rf}{r}-\frac{f}{r^2})\right)\underline{\omega}P_n\\
                                &=&((\Delta_zf+(2n+\mu-1)D^r(1)\{f\})\underline{\omega}P_n.
\end{eqnarray*}
Summarizing we have that the lemma is true in the case $m=1$. Assume
that our formulae hold for a positive integer $m$, we have to show
them for $m+1$.

We thus get
\begin{eqnarray*}
\Delta^{m+1}(fP_n)&=&\sum_{j=0}^md_{n,\mu}(j)\left(\begin{array}
{cc} m \\ j
\end{array}\right)\Delta\left(D_r(j)\{\Delta_z^{m-j}f\}P_n\right)\\
                    &=&\sum_{j=0}^md_{n,\mu}(j)\left(\begin{array}
{cc} m \\ j
\end{array}\right)(\partial^2_{x_0}D_r(j)\{\Delta_z^{m-j}f\}+\partial^2_rD_r(j)\{\Delta_z^{m-j}f\}\\
                    & &+\frac{2n+\mu-1}{r}\partial_rD_{r}(j)\{\Delta_z^{m-j}f\})P_n\\
                    &=&\sum_{j=0}^md_{n,\mu}(j)\left(\begin{array}
{cc} m \\ j
\end{array}\right)(D_r(j)\{\partial^2_{x_0}\Delta_z^{m-j}f+\partial^2_r\Delta_z^{m-j}f\}\\
                    & &-2jD_r(j+1)\{\Delta_z^{m-j}f\}+(2n+\mu-1)D_{r}(j+1)\Delta_z^{m-j}f)P_n\\
                    &=&\left(\sum_{j=0}^md_{n,\mu}(j)\left(\begin{array}
{cc} m \\ j
\end{array}\right)D_r(j)\{\Delta_z^{m+1-j}f\}\right)P_n\\
                    & &+\left(\sum_{j=0}^md_{n,\mu}(j)\left(\begin{array}
{cc} m \\ j
\end{array}\right)(2n+\mu-(2j+1))D_r(j+1)\{\Delta_z^{m-j}f\}\right)P_n\\
                    &=&\left(\sum_{j=0}^{m+1}d_{n,\mu}(j)\left(\begin{array}
{cc} m+1 \\ j
\end{array}\right)D_r(j)\{\Delta_z^{m+1-j}f\}\right)P_n.
\end{eqnarray*}
In the last step we used the fact that $\left(\begin{array} {cc} m
\\ j
\end{array}\right)+\left(\begin{array}
{cc} m \\ j-1
\end{array}\right)=\left(\begin{array}
{cc} m+1 \\ j
\end{array}\right)$. This establishes the
first formula. The other one may be proved in a similar way. $\qquad
\blacksquare$

Now we arrive at the following higher order version of Fueter's
Theorem in Dunkl case.

\begin{theorem}
Let $f(z)=u(t,s)+\i v(t,s)$, $z=t+\i s$, be a complex-valued
function satisfying equation (\ref{Eqn:1}) in some open set
$\Omega\subset\C^+=\{z\in\C: s>0\}$. If $\mu$ is odd, then the
following function
\begin{eqnarray*}
\Delta^{m+n+\frac{\mu-1}{2}}\left(\left(u(x_0,|\underline{x}|)+\frac{\underline{x}}{|\underline{x}|}v(x_0,|\underline{x}|)\right)P_n(\underline{x})\right)
\end{eqnarray*}
is Dunkl-monogenic in $\vec{\Omega}=\{x\in\R_1^d:
(x_0,|\underline{x}|)\in\Omega\}$.
\end{theorem}

\proof From Lemma 3.2, we get that
\begin{eqnarray*}
&
&\Delta^{m+n+\frac{\mu-1}{2}}(u(x_0,|\underline{x}|)P_n(\underline{x}))\\
&=&\left(\sum_{j=0}^{m+n+\frac{\mu-1}{2}}d_{n,\mu}(j)\left(\begin{array}
{cc} m+n+\frac{\mu-1}{2} \\ j
\end{array}\right)D_r(j)\{\Delta_z^{m+n+\frac{\mu-1}{2}-j}u(x_0,r)\}\right)P_n(\underline{x})
\end{eqnarray*}
and
\begin{eqnarray*}
&
&\Delta^{m+n+\frac{\mu-1}{2}}(v(x_0,|\underline{x}|)\frac{\underline{x}}{|\underline{x}|}P_n(\underline{x}))\\
&=&\left(\sum_{j=0}^{m+n+\frac{\mu-1}{2}}d_{n,\mu}(j)\left(\begin{array}
{cc} m+n+\frac{\mu-1}{2} \\ j
\end{array}\right)D^r(j)\{\Delta_z^{m+n+\frac{\mu-1}{2}-j}v(x_0,r)\}\right)\frac{\underline{x}}{|\underline{x}|}P_n(\underline{x}).
\end{eqnarray*}
Obviously, the first $n+(\mu-1)/2$ terms in the above equalities
vanish since by hypothesis $f$ satisfies equation (\ref{Eqn:1}).
Furthermore, note that $2n+\mu-(2j-1)<0$ for $j\geq n+(\mu+1)/2$ and
therefore $d_{n,\mu}(j)=0$ for $j\geq n+(\mu+1)/2$.

Summarizing we obtain that
\begin{eqnarray*}
&
&\Delta^{m+n+\frac{\mu-1}{2}}\left(\left(u(x_0,|\underline{x}|)+\frac{\underline{x}}{|\underline{x}|}v(x_0,|\underline{x}|)\right)P_n(\underline{x})\right)\\
&=&(2n+\mu-1)!!\left(\begin{array} {cc} m+n+\frac{\mu-1}{2} \\
n+\frac{\mu-1}{2}
\end{array}\right)(A(x_0,r)+\frac{\underline{x}}{|\underline{x}|}B(x_0,r))P_n(\underline{x}),
\end{eqnarray*}
where
\begin{eqnarray*}
& &A=D_r\left(n+\frac{\mu-1}{2}\right)\{\Delta_z^mu\},\\
& &B=D^r\left(n+\frac{\mu-1}{2}\right)\{\Delta_z^mv\}.
\end{eqnarray*}
The task is now to show that $A$ and $B$ satisfy the following
Vekua-type system in Dunkl case(\cite{FCK3})
\begin{eqnarray*}
\left\{
       \begin{array}{ll}
           \partial_{x_0}A-\partial_rB=\frac{2n+\mu-1}{r}B,\\
        \ \\
           \partial_{x_0}B+\partial_rA=0.\\
        \end{array}
        \right.
\end{eqnarray*}
To this end it will be necessary to use the assumptions on $u$ and
$v$
\begin{eqnarray*}
\left\{
       \begin{array}{ll}
           \partial_t\Delta_z^mu-\partial_s\Delta_z^mv=0,\\
        \ \\
           \partial_t\Delta_z^mv+\partial_s\Delta_z^mu=0,\\
        \end{array}
        \right.
\end{eqnarray*}
and statements ${\rm (iii)}$ and ${\rm (iv)}$ of Lemma 3.1.

Indeed, we have
\begin{eqnarray*}
\partial_{x_0}A-\partial_rB&=&D_r\left(n+\frac{\mu-1}{2}\right)\{\partial_{x_0}\Delta_z^mu\}-\partial_rD^r\left(n+\frac{\mu-1}{2}\right)\{\Delta_z^mv\}\\
                           &=&D_r\left(n+\frac{\mu-1}{2}\right)\{\partial_r\Delta_z^mv\}-\partial_rD^r\left(n+\frac{\mu-1}{2}\right)\{\Delta_z^mv\}\\
                           &=&\frac{2n+\mu-1}{r}D^r\left(n+\frac{\mu-1}{2}\right)\{\Delta_z^mv\}\\
                           &=&\frac{2n+\mu-1}{r}B
\end{eqnarray*}
and
\begin{eqnarray*}
\partial_{x_0}B+\partial_rA&=&D^r\left(n+\frac{\mu-1}{2}\right)\{\partial_{x_0}\Delta_z^mv\}+\partial_rD_r\left(n+\frac{\mu-1}{2}\right)\{\Delta_z^mu\}\\
                           &=&D^r\left(n+\frac{\mu-1}{2}\right)\{\partial_{x_0}\Delta_z^mv\}+D^r\left(n+\frac{\mu-1}{2}\right)\{\partial_r\Delta_z^mu\}\\
                           &=&D^r\left(n+\frac{\mu-1}{2}\right)\{\partial_{x_0}\Delta_z^mv+\partial_r\Delta_z^mu\}\\
                           &=&0,
\end{eqnarray*}
which completes the proof. $\qquad \blacksquare$

\section{Fueter's Theorem with an extra Dunkl-monogenic factor $P_n(x_0,\underline{x})$}

We begin this section with two basic theorems of Dunkl-Clifford
analysis. The proof of Theorem 4.1 is straightforward from the one
of CK-Extension Theorem in classical Clifford analysis. Theorem 4.2
is from~\cite{Ren} or~\cite{BCK}.
\begin{theorem}
Any analytic function $f(\underline{x})$ in $\R^d$ has a unique
Dunkl-monogenic extension $CK[g]$ to $\R_1^d$, which is given by
$$CK[g(\underline{x})](x)=\sum_{j=0}^{\infty}\frac{(-x_0)^j}{j!}\underline{D}^jg(\underline{x}).$$
\end{theorem}

Let ${\cal P}(n)$, $n\in\N$, denote the set of all $\R_{0,d}-$valued
homogeneous polynomials of degree $n$ in $\R^d$, which contain the
important subspace ${\cal M}(n)$ introduced in Section 2. Then there
holds the following Almasi-Fischer Decomposition Theorem.
\begin{theorem}
Let $n\in\N$. Then $${\cal
P}(n)=\bigoplus_{k=0}^n\underline{x}^k{\cal M}(n-k).$$
\end{theorem}

Now we are ready to prove the following Fueter's Theorem in
Dunkl-Clifford analysis with an extra Dunkl-monogenic factor
$P_n(x_0,\underline{x})$.
\begin{theorem}
Let $f(z)=u(t,s)+\i v(t,s)$ be a complex-valued holomorphic function
in some open set $\Omega\subset\C^+$ and assume that
$P_n(x_0,\underline{x})$ is a homogeneous Dunkl-monogenic polynomial
of degree $n$ in $\R_1^d$. If $\mu$ is odd, then the following
 function
\begin{eqnarray*}
\Delta^{n+\frac{\mu-1}{2}}\left(\left(u(x_0,|\underline{x}|)+\frac{\underline{x}}{|\underline{x}|}v(x_0,|\underline{x}|)\right)P_n(x_0,\underline{x})\right)
\end{eqnarray*}
is Dunkl-monogenic in $\vec{\Omega}=\{x\in\R_1^d:
(x_0,|\underline{x}|)\in\Omega\}$.
\end{theorem}

\proof It is obvious from Theorem 4.1 that
$P_n(x_0,\underline{x})=CK[P_n(0,\underline{x})](x)$. By Theorem 4.2
there exist unique $P_{n-k}\in{\cal M}(n-k)$ such that
$$P_n(x_0,\underline{x})=\sum_{k=0}^nCK[\underline{x}^kP_{n-k}(\underline{x})](x).$$
Thus, it suffices to show that the Dunkl-monogenicity of function
\begin{eqnarray*}
\Delta^{n+\frac{\mu-1}{2}}\left(\left(u(x_0,|\underline{x}|)+\frac{\underline{x}}{|\underline{x}|}v(x_0,|\underline{x}|)\right)CK[\underline{x}^kP_{n-k}(\underline{x})](x)\right),
\qquad k=0,\cdots,n.
\end{eqnarray*}
Since
\begin{eqnarray*}
CK[\underline{x}^kP_{n-k}(\underline{x})](x)=\sum_{j=0}^n\frac{(-x_0)^j}{j!}\underline{D}^j(\underline{x}^kP_{n-k}(\underline{x}))
\end{eqnarray*}
and by means of Lemma 2.1, we can conclude that
$CK[\underline{x}^kP_{n-k}(\underline{x})](x)$ is of the form
\begin{eqnarray*}
CK[\underline{x}^kP_{n-k}(\underline{x})](x)=\left(\sum_{j=0}^kc_jx_0^j\underline{x}^{k-j}\right)P_{n-j}(\underline{x}),
\ c_j\in\R.
\end{eqnarray*}
Therefore,
\begin{eqnarray*}
CK[\underline{x}^kP_{n-k}(\underline{x})](x)=\left(U(x_0,r)+\frac{\underline{x}}{|\underline{x}|}V(x_0,r)\right)P_{n-k}(\underline{x}),
\end{eqnarray*}
where $U$ and $V$ are real-valued homogeneous polynomials of degree
$k$ in the variables $x_0$ and $r$. So, its corresponding
complex-valued function $g(z)=U(t,s)+\i V(t,s)$ obviously satisfies
\begin{eqnarray*}
\partial_{\bar{z}}^{k+1}g(z)=0,\ z=t+\i s\in\C.
\end{eqnarray*}
Whence by the assumption of $f$,
\begin{eqnarray*}
\partial_{\bar{z}}^{k+1}(f(z)g(z))=0,\ z\in\Omega,
\end{eqnarray*}
i.e. $f(z)g(z)$ is $(n+1)-$holomorphic in $\Omega$. It then follows
that
\begin{eqnarray*}
\partial_{\bar{z}}\Delta_z^{k}(f(z)g(z))=0,\ z\in\Omega.
\end{eqnarray*}
The proof of Theorem 4.3 now follows by using Theorem 3.1. $\qquad
\blacksquare$

\begin{remark}
Comparing with the classical case(see~\cite{PS2010}), the crucial
part in our treatment of Theorem 4.3 is the replacement of the
classical Euclidean dimension $d$ by the Dunkl-dimension $\mu$.
\end{remark}

\noindent{\large \bf  Acknowledgements}

\bigskip

\noindent Part of the work in this paper was done while the authors
were visiting at University of Aveiro during the academic year
2009-10. The authors were (partially) supported by {\it CIDMA -
Centro de Investiga\c c\~ao e Desenvolvimento em Matem\'atica e
Aplica\c c\~oes} of the University of Aveiro. The second author is
the recipient of a grant from {\it Funda\c{c}\~{a}o para Ci\^{e}ncia
e a Tecnologia (Portugal)} with grant No.: SFRH/BPD/41730/2007.

 \end{document}